\documentclass[11pt,twoside]{article}

\usepackage{latexsym}
\usepackage{amsmath}
\usepackage{amsthm}
\usepackage{amssymb}
\usepackage{vmargin}
\usepackage{amscd}
\usepackage{stmaryrd}
\usepackage{euscript}
\usepackage{mathrsfs}
\usepackage{amscd}
\usepackage[all]{xy}
\usepackage{xr}
\DeclareMathAlphabet{\mathpzc}{OT1}{pzc}{m}{it}

\setmargins{32mm}{20mm}{14.6cm}{22cm}{1cm}{1cm}{1cm}{1cm}

\def\ra{\rightarrow}
\def\la{\leftarrow}
\def\lra{\longrightarrow}
\def\lla{\longleftarrow}
\def\id{\mathrm{id}}
\def\s|{\,|\,}
\def\ten{\otimes}
\def\vareps{\varepsilon}
\def\eps{\epsilon}

\def\CC{\mathrm{C}}

\def\Ru{\mathrm{R_u}}

\def\H{\mathrm{H}}
\def\z{\mathrm{Z}}

\def\Z{\mathbb{Z}}

\def\R{\mathbb{R}}
\def\Cx{\mathbb{C}}

\def\vv{\mathbb{V}}
\def\ww{\mathbb{W}}
\def\Bu{\mathbb{B}}

\def\bB{\mathbb{B}}

\def\bG{\mathbb{G}}

\def\bO{\mathbb{O}}

\def\C{\mathcal{C}}

\def\cL{\mathcal{L}}

\def\cN{\mathcal{N}}

\def\sA{\mathscr{A}}

\def\Def{\mathfrak{Def}}

\def\fB{\mathfrak{B}}

\def\fR{\mathfrak{R}}

\def\L{\Lambda}

\def\m{\mathfrak{m}}

\def\g{\mathfrak{g}}
\def\ga{\mathpzc{g}}

\def\fh{\mathfrak{h}}

\def\fu{\mathfrak{u}}

\def\Hom{\mathrm{Hom}}

\def\Aut{\mathrm{Aut}}

\def\Iso{\mathrm{Iso}}

\def\coker{\mathrm{coker\,}}

\def\Spec{\mathrm{Spec}\,}

\def\Set{\mathrm{Set}}

\def\Grp{\mathrm{Grp}}

\def\Grpd{\mathrm{Grpd}}

\def\ad{\mathrm{ad}}

\def\<{\langle}
\def\>{\rangle}
\def\Lim{\varprojlim}
\def\into{\hookrightarrow}
\def\onto{\twoheadrightarrow}

\def\xra{\xrightarrow}

\def\alg{\mathrm{alg}}

\def\by{\times}

\def\mc{\mathrm{MC}}

\def\ddef{\mathrm{Def}}

\def\mca{\mathrm{MC}_{A}}

\def\ga{\mathrm{G}_{A}}

\def\defa{\mathrm{Def}_{A}}

\def\DGVect{\mathrm{DGVect}}
\def\Vect{\mathrm{Vect}}
\def\Rep{\mathrm{Rep}}

\def\SL{\mathrm{SL}}
\def\GL{\mathrm{GL}}

\def\ind{\mathrm{ind}}
\def\pro{\mathrm{pro}}

\def\iff{\Leftrightarrow}

\def\pd{\partial}
\def\dc{d^{\mathrm{c}}}
\def\half{\frac{1}{2}}

\def\gr{\mathrm{gr}}

\def\redpi{\varpi_1^{\mathrm{red}}(X,x)}
\def\algpi{\varpi_1(X,x)}
\def\mypi{R_u(\varpi_1(X,x))}

\def\alg{\mathrm{alg}}
\def\red{\mathrm{red}}

\newtheorem{theorem}{Theorem}[section]
\newtheorem{proposition}[theorem]{Proposition}
\newtheorem{corollary}[theorem]{Corollary}
\newtheorem{lemma}[theorem]{Lemma}
\newtheorem{theorem*}{Theorem}
\newtheorem{proposition*}[theorem*]{Proposition}
\newtheorem{corollary*}[theorem*]{Corollary}
\newtheorem{lemma*}[theorem*]{Lemma}

\theoremstyle{definition}
\newtheorem{definition}[theorem]{Definition}

\newtheorem{definition*}[theorem*]{Definition}

\theoremstyle{remark}
\newtheorem{example}[theorem]{Example}

\newtheorem{remark}[theorem]{Remark}
\newtheorem{remarks}[theorem]{Remarks}

\newtheorem{example*}[theorem*]{Example}
\newtheorem{examples*}[theorem*]{Examples}
\newtheorem{remark*}[theorem*]{Remark}
\newtheorem{remarks*}[theorem*]{Remarks}
\newtheorem{exercise*}[theorem*]{Exercise}

\begin{document}
\title{The pro-unipotent radical of the  pro-algebraic fundamental group of a compact K\"ahler manifold}
\author{J. P. Pridham\thanks{The author is supported by Trinity College, Cambridge}}

\maketitle

\tableofcontents

\section*{Introduction}
\addcontentsline{toc}{section}{Introduction}

For $X$  a compact connected K\"ahler manifold, consider the real pro-algebraic completion $\algpi$ of the  fundamental group $\pi_1(X,x)$. In \cite{Simpson}, Simpson defined a Hodge structure on the complex pro-algebraic fundamental group $\algpi_{\Cx}$, in the form of a discrete $\Cx^*$-action. The Levi decomposition for pro-algebraic groups allows us to write
$$
\algpi \cong \mypi\rtimes \redpi,
$$
where  $\mypi$ is the pro-unipotent radical of $\algpi$ and $\redpi$ is the reductive quotient of $\algpi$. This decomposition is  unique up to conjugation by $\mypi$. 
By studying the Hodge structure  on $\redpi$, Simpson established restrictions on its possible group structures. 

The purpose of this paper is to use Hodge theory to show that  $\mypi$ is quadratically presented as a pro-unipotent group, in the sense that its Lie algebra can be defined by equations of bracket length two.   This  generalises both Goldman and Millson's result on deforming reductive representations of the fundamental group (\cite{GM}), and Deligne et al.'s result on the de Rham fundamental group $\pi_1(X,x)\ten \R$  (\cite{DGMS}). The idea behind this paper is that in both of these cases, we are taking a reductive representation
$$
\rho_0:\pi_1(X,x) \to G,
$$
and considering deformations
$$
\rho: \pi_1(X,x) \to U \rtimes G
$$
of $\rho_0$, for $U$ unipotent.    

Effectively, \cite{GM} considers only $U=\exp(\mathrm{Lie}(G) \ten \m_A)$, for  $\m_A$ a maximal ideal of an Artinian local $\R$-algebra, while \cite{DGMS} considers only the case when $G=1$. Since taking $U=\mypi$ pro-represents this functor when $G=\redpi$, the quadratic presentation for $U$ gives quadratic presentations both for the hull of \cite{GM} and for the Lie algebra of \cite{DGMS}.  This also generalises Hain's results (\cite{malcev}) on relative Malcev completions of variations of Hodge structure, since here we are taking relative Malcev completions of arbitrary reductive representations.

Section \ref{review} summarises standard definitions and properties of pro-algebraic groups which are used throughout the rest of the paper.

Section \ref{nilplie} develops a theory of deformations over nilpotent Lie algebras with $G$-actions. This can be thought of as a generalisation of the theory introduced in \cite{pi1}, which corresponds to the case $G=1$. The essential philosophy is that all the concepts for deformations over Artinian rings, developed by Schlessinger in \cite{Sch} and later authors, can be translated to this context.

Section \ref{tdga} introduces the notion of twisted differential graded algebras, which are analogous to the DGAs used in \cite{DGMS} to characterise the real homotopy type. They are equivalent to the $G$-equivariant DGAs used in \cite{KTP} to study the schematic homotopy type.

Section \ref{proalg} contains various technical lemmas about pro-algebraic groups. 

In Section \ref{torsor}, the twisted DGA arising from  $\CC^{\infty}$-sections is defined.
It is shown that this can be used to recover  $\mypi$. 

Section \ref{hodge} uses Hodge theory to prove that, for a compact K\"ahler manifold, this twisted DGA is formal, i.e. quasi-isomorphic to its cohomology DGA. This can be thought of as equivalent to formality of the real schematic homotopy type.   In consequence, $\mypi$ is quadratically presented --- a direct analogue of the demonstration in \cite{DGMS} that formality of the real homotopy type forces the de Rham fundamental group to be quadratically presented. 
This implies that if $\pi_1(X,x)$ is of the form $\Delta \rtimes \Lambda$, with $\Lambda$ acting reductively on $\Delta \ten \R$, then $\Delta \ten \R$ must be quadratically presented. 

I would like to thank Bertrand To\"en for suggesting the connection between this work and the schematic homotopy type.

\section{Review of pro-algebraic groups}\label{review}
In this section, we recall some definitions and standard properties of pro-algebraic groups, most of which can be found in \cite{tannaka} II\S 2 and \cite{Simpson}. Fix a field $k$ of characteristic zero. 

\begin{definition}
An algebraic group over $k$ is defined to be an affine group scheme $G$ of finite type over $k$. These all arise as Zariski-closed subgroups of general linear groups $\GL_n(k)$. A pro-algebraic group is a filtered inverse limit of algebraic groups, or equivalently an arbitrary affine group scheme over $k$.   
\end{definition}

\begin{definition}\label{O}
Given a pro-algebraic group $G$, let $O(G)$ denote global sections of the structure sheaf of $G$, so that $G=\Spec O(G)$. This is a sum of  finite-dimensional $G\by G$-representations, the actions corresponding to right and left translation. The group structure on $G$ corresponds to a comultiplication $\Delta:O(G)\to O(G)\ten O(G)$, coidentity $\vareps:O(G)\to k$, and coinverse $S:O(G)\to O(G)$, satisfying coassociativity, coidentity and coinverse axioms.
\end{definition}

\begin{lemma}\label{dual} If $G$ is a pro-algebraic group, and we regard $O(G)$ as a sum of finite-dimensional $G$-representations via the left action, then for any finite-dimensional $G$-representation $V$, 
$$
V\cong (V\ten O(G))^G:= \{a \in V\ten O(G)\,|\, (g\ten \id)a=(\id \ten g)a,\,  \forall g \in G\},
$$
with the $G$-action on the latter coming from the right action on $O(G)$.
\end{lemma}
\begin{proof}
This follows immediately from  \cite{tannaka} II Proposition 2.2, which states that $G$-representations correspond to $O(G)$-comodules. Under this correspondence, $v \in V$ is associated to  the function $g \mapsto g\cdot v$. 
\end{proof}

\begin{definition}
An algebraic group $G$ is said to be unipotent if the coproduct $\Delta:O(G)\to O(G)\ten O(G)$ is counipotent. Unipotent algebraic groups all arise as Zariski-closed subgroups of the groups of upper triangular matrices:
$$
\begin{pmatrix} 1 & * & *\\0 & \ddots & * \\ 0 & 0 & 1 \end{pmatrix}. 
$$

A pro-algebraic group is said to be pro-unipotent if it is an inverse limit of unipotent algebraic groups. This is equivalent to saying that $\Delta$ is ind-counipotent.
\end{definition}

\begin{lemma}
There is a one-to-one correspondence between unipotent algebraic groups over $k$, and finite-dimensional nilpotent Lie algebras over $k$.
\end{lemma}
\begin{proof}
The Lie algebra $\fu$ of  any unipotent algebraic group $U$ is necessarily finite-dimensional and nilpotent. Conversely, if $\fu$ is any finite-dimensional nilpotent Lie algebra, we define a unipotent algebraic group $\exp(\fu)$ by
$$
\exp(\fu)(A):= \exp(\fu\ten A).
$$
Here, for any Lie algebra $\g$, the group $\exp(\g)$ has underlying set $\g$ and multiplication given by the Campbell-Baker-Hausdorff formula
$$
g\cdot h:= g +h +\half [g,h] +\ldots,
$$
which in this case is a finite sum, by nilpotence, so the group $\exp(\fu)$ is indeed algebraic. To see that these functors are inverse is most easily done by considering groups of upper triangular matrices.
\end{proof}

\begin{definition}
The pro-unipotent radical $\Ru(G)$ of a pro-algebraic group $G$ is defined to the maximal pro-unipotent normal closed subgroup of $G$. A pro-algebraic group $G$ is said to be reductive if $\Ru(G)=1$, and for an arbitrary pro-algebraic group $G$, the reductive quotient of $G$ is defined by $G^{\red}:=G/\Ru(G)$.
\end{definition}

\begin{theorem}[Levi decomposition]\label{radical} For any pro-algebraic group $G$, there is 
a decomposition
$$
G \cong\Ru(G) \rtimes G^{\red},
$$ 
unique up to conjugation by the pro-unipotent radical  $\Ru(G)$.  
\begin{proof}
This is the Levi decomposition for pro-algebraic groups in characteristic zero, proved in \cite{Levi}, which states that for every pro-algebraic group $G$, the surjection $G \to G^{\red}$ has a section, unique up to conjugation by $\Ru(G)$, inducing an isomorphism $G\cong \Ru(G)\rtimes G^{\red}$. 
\end{proof}
\end{theorem}

\begin{theorem}[Tannakian duality]
A pro-algebraic group $G$ over $k$  can be recovered from its (tensor) category of finite-dimensional $k$-representations. Representations of $G^{\red}$ correspond to the subcategory of semisimple representations.
\end{theorem}
\begin{proof}
The first part is \cite{tannaka} Theorem II.2.11. The second part is just the observation that, in characteristic zero, ``reductive'' and ``linearly reductive'' are equivalent.
\end{proof}

\begin{definition}
Given a discrete group $\Gamma$, define the pro-algebraic completion $\Gamma^{\alg}$ of $\Gamma$ to be the pro-algebraic group $G$ universal among group homomorphisms
$
\Gamma \to G(k).
$
In other words $\Gamma^{\alg}$ pro-represents the functor which sends an algebraic group $G$ to the set of representations  $\Gamma \to G(k)$. Under Tannakian duality, the finite-dimensional linear $k$-representations of $\Gamma^{\alg}$ are just the finite-dimensional linear $k$-representations of $\Gamma$.

The reductive quotient of $\Gamma^{\alg}$ is denoted $\Gamma^{\red}$, and is  universal among Zariski-dense group homomorphisms
$
\Gamma \to G(k),
$
with $G$ reductive. In other words $\Gamma^{\red}$ pro-represents the functor which sends an algebraic group $G$ to the set of reductive representations  $\Gamma \to G(k)$. Under Tannakian duality, the finite-dimensional linear $k$-representations of $\Gamma^{\red}$ are just the semisimple finite-dimensional linear $k$-representations of $\Gamma$.

The pro-unipotent (or Malcev) completion $\Gamma \ten k$ of $\Gamma$ is  universal among group homomorphisms
$
\Gamma \to G(k),
$
with $G$ pro-unipotent. In other words $\Gamma\ten k$ pro-represents the functor which sends an algebraic group $G$ to the set of unipotent representations $\Gamma \to G(k)$. Under Tannakian duality, the finite-dimensional linear $k$-representations of $\Gamma\ten k $ are just the unipotent  finite-dimensional linear $k$-representations of $\Gamma$.
\end{definition}

\begin{definition} Given a representation $\rho:\Gamma \to G(k)$, for some pro-algebraic group $G$,  the Malcev completion of $\Gamma$ relative to $\rho$ is defined in \cite{malcev} to represent the functor which sends a  pro-unipotent extension $H \to G$ to the set of representations
$$
\Gamma \to H(k)
$$
lifting $\rho$. Thus $\Gamma\ten k$ is the Malcev completion of $\Gamma$ relative to the trivial representation, and $\Gamma^{\alg}$ is the Malcev completion of $\Gamma$ relative to the canonical representation $\Gamma \to \Gamma^{\red}$.
\end{definition}

\section{Functors on nilpotent Lie algebras with $G$-actions}\label{nilplie}

This section extends the ideas of \cite{Sch} to a slightly different context.

Fix a field $k$ of characteristic zero. Take a  pro-algebraic group $G$ over $k$, and let $\Rep(G)$ be the category of finite-dimensional  representations of $G$ over $k$. If $G$ is reductive, then every such representation will be decomposable into irreducibles, so $\Hom$ will be an exact functor on this category. Consider the category $\widehat{\Rep}(G):=\pro(\Rep(G))$, whose objects are filtered inverse systems $\{V_{\alpha}\}_{\alpha \in I}$, with morphisms given by
$$
\Hom_{\pro(\Rep(G))}(\{V_{\alpha}\}, \{W_{\beta}\})= \lim_{\substack{\lla \\ \beta}} \lim _{\substack{\lra \\ \alpha}}\Hom_{\Rep(G)}(V_{\alpha},W_{\beta}).
$$

Given a set $\{V_i\}_{i \in I}$, with $V_i \in \Rep(G)$, we make the vector space $\prod_{i\in I} V_i$ an object of $\widehat{\Rep}(G)$ via the formula
$$
\prod_{i\in I} V_i= \lim_{\substack{\lla \\ J \subset I \text{ finite}}} \prod_{j\in J} V_j.
$$

\begin{lemma} If $G$ is reductive, then every object of $\widehat{\Rep}(G)$ can be expressed as a product of irreducible finite-dimensional $G$-representations.
\begin{proof}
Since $\Rep(G)$ is an Artinian category, i.e. it satisfies the descending chain condition for sub-objects, we may use \cite{descent} to observe that $\widehat{\Rep}(G)$ is isomorphic to the category of left-exact set-valued functors on $\Rep(G)$.

Take $W \in \widehat{\Rep}(G)$, and let $\{V_s: s \in S\}$ be a set of representatives for isomorphism classes of irreducible representations in $\Rep(G)$. Now, $\Hom_{\widehat{\Rep}(G)}(W,V_s)$ has the natural structure of a vector space over $k$. Choose a basis $t_i: i \in I_s$ for this vector space, and let 
$$
U:= \prod_{s \in S} V_s^{I_s}.
$$
There is then a natural isomorphism between $ \Hom_{\widehat{\Rep}(G)}(W,V_s)$ and $\Hom_{\widehat{\Rep}(G)}(U,V_s)$ for all $s \in S$, so the left-exact functors defined on $\Rep(G)$ by $U$ and $W$ must be isomorphic, and therefore $U \cong W$.
\end{proof}
\end{lemma}

\begin{definition}
For any pro-algebraic group $G$, define $\cN(G)$ to be the category whose objects are  pairs $(\fu, \rho)$, where $\fu$ is a  finite-dimensional nilpotent Lie algebras over $k$, and $\rho: G \to \Aut(\fu)$ is a representation to the group of Lie algebra automorphisms of $\fu$. A morphism $\theta$ from $(\fu,\rho)$ to $(\fu', \rho')$ is a morphism  $\theta: \fu \to \fu'$ of Lie algebras such that $\theta\circ\rho=\rho'$. Observe that $\cN(G)$ is an Artinian category,  and write $\widehat{\cN}(G)$ for the category $\pro(\cN(G))$.
\end{definition}

Given $\cL \in \widehat{\cN}(G)$, let $\cN(G)_{\cL}$ be the category of pairs $(N \in \cN(G), \cL \xra{\phi} N)$, and $\widehat{\cN}(G)_{\cL}$ the category of pairs $(N \in \widehat{\cN}(G), \cL \xra{\phi} N)$. We will almost always consider the case $\cL=0$ (note that $\cN(G)_0=\cN(G)$), although a few technical lemmas (Propositions \ref{nMan2} and \ref{nMan3} and Corollary \ref{nkeyhgs}) require the full generality of $\cN(G)_{\cL}$. In $\widehat{\cN}(G)_{\cL}$, $\cL$ is the initial object, and $0$ the final object.

\begin{definition} For $N = \Lim_{\alpha} N_{\alpha} \in \widehat{\cN}(G)_{\cL}$, define the tangent space functor 
\begin{eqnarray*}
t_{N/\cL}: \Rep(G) &\to& \Vect\\
V &\mapsto& \Hom_{\widehat{\cN}(G)_{\cL}}(N, V\eps), 
\end{eqnarray*}
    where $V\eps \in \cN(G)$ is the representation  $V$  regarded as an abelian Lie algebra, i.e. \mbox{$[V\eps,V\eps]=0$}, with structure morphism $\cL \xra{0} V\eps$. This is clearly a vector space over $k$, and 
$$
t_{N/\cL}(V)=\lim_{\substack{\lra \\ \alpha}}\Hom_{\Rep(G)}( N_{\alpha}/\langle[N_{\alpha},N_{\alpha}] +\cL\rangle , V).  
$$
We define the cotangent space
$$
t_{N/\cL}^*:=N/\langle[N,N] +\cL\rangle = \lim_{\substack{\lra \\ \alpha}} N_{\alpha}/\langle[N_{\alpha},N_{\alpha}] +\cL\rangle \in \widehat{\Rep}(G).
$$ 
\end{definition}

\begin{definition}\label{nfreelie} Given $V$ in $\widehat{\Rep}(G)$, denote the free pro-nilpotent Lie algebra on generators $V$ by $L(V)$. This has a natural continuous $G$-action, so is in $\widehat{\cN}(G)$. Equivalently, we may use \cite{descent} to define $L(V)$ as the object of $\widehat{\cN}(G)$ pro-representing the functor $N \mapsto \Hom_{\widehat{\Rep}(G)}(V,N)$.
\end{definition}

\begin{definition} Given a Lie algebra $L$, we define the lower central series of ideals $\Gamma_n(L)$  inductively by
$$
\Gamma_1(L)=L, \quad \Gamma_{n+1}(L)=[L,\Gamma_n(L)],
$$
and define the associated graded algebra $\gr L$ by $\gr_n L= \Gamma_n L/\Gamma_{n+1} L$ 

For $L \in \cN(G)_{\cL}$, define the ideals $\Gamma^{\cL}_n(L)$ by
$$
\Gamma_1(L)=L, \quad \Gamma_2(L)=[L,l]+\langle \cL \rangle \quad \Gamma_{n+1}(L)=[L,\Gamma_n(L)],
$$
and define the associated graded algebra $\gr^{\cL} L$ by $\gr^{\cL}_n L= \Gamma^{\cL}_n L/\Gamma^{\cL}_{n+1} L$ 
\end{definition}

\begin{definition}
Given Lie algebras $N,M \in \widehat{\cN}(G)$, let $N*M$ be the completed free Lie algebra product, i.e. the completion with respect to the commutator filtration of the free product of the pro- Lie algebras $N$ and $M$. 
 Note that $*$ is sum in $\widehat{\cN}(G)$ --- the analogue in the category $\pro(\C_{\Lambda})$ of pro-Artinian $\Lambda$-algebras is $\hat{\ten}$.

Given  Lie algebras $N,M \in \widehat{\cN}(G)_{\cL}$, define the free fibre product $N*_{\cL}M$ similarly. Equivalently, we can use \cite{descent} to define this element of $ \widehat{\cN}(G)$, since all finite colimits must exist in a pro-Artinian category.
\end{definition}

\begin{definition}
Define the free Lie algebra $L_{\cL}(V):=\cL * L(V)$
\end{definition}

We will consider only those  functors $F$ on $\cN(G)_{\cL}$ which satisfy 
\begin{enumerate}
\item[(H0)]$F(0)=\bullet$, the one-point set.
\end{enumerate}

We adapt the following definitions and results from \cite{Sch} (with identical proofs):

\begin{definition} For $p:N \ra M$ in $\widehat{\cN}(G)_{\cL}$ surjective, $p$ is a semi-small extension if  $[N,\ker p]=(0)$. If the pro-$G$-representation $\ker p$ is an absolutely irreducible $G$-representation, then we say that $p$ is a small extension. Note that any surjection in $\cN(G)_{\cL}$ can be be factorised as a composition of small extensions.
\end{definition}

For $F:\cN(G)_{\cL} \ra \Set$, define $\hat{F}:\widehat{\cN}(G)_{\cL} \ra \Set$ by 
$$
\hat{F}(\lim_{\substack{\longleftarrow \\ \alpha}}  L_{\alpha})= \lim_{\substack{\longleftarrow \\ \alpha}} F(L_{\alpha}).
$$

Note that \mbox{$\hat{F}(L) \xrightarrow{\sim} \Hom(h_L,F)$,} where 
\begin{eqnarray*}
h_L:\cN(G)_{\cL} &\to& \Set;\\ 
N &\mapsto& \Hom(L,N).
\end{eqnarray*}

\begin{definition}
We will say a functor $F:\cN(G)_{\cL} \ra \Set$ is pro-representable if it is isomorphic to $h_L$, for some $L \in \widehat{\cN}(G)_{\cL}$. By the above remark, this isomorphism is determined by an element $\xi \in \hat{F}(L)$. We say the pro-couple $(L,\xi)$ pro-represents $F$.
\end{definition}

\begin{remark}
This definition is not the strict analogue of that appearing in \cite{Sch}, which had additional hypotheses on the finiteness of tangent spaces. This terminology coincides with that of \cite{descent}.
\end{remark}

\begin{definition} A natural transformation $\phi:F \ra E$ in $[\cN(G)_{\cL},\Set]$ is called:
\begin{enumerate}
\item unramified if $\phi:F(V\eps) \ra E(V\eps)$ is injective for all irreducible $G$-representations $V$.

\item smooth if for every surjection $M \twoheadrightarrow N$ in $\cN(G)_{\cL}$, the canonical map  \mbox{$F(M) \to E(M)\by_{E(N)}F(N)$} is surjective.
\item \'etale if it is smooth and unramified.
\end{enumerate}
\end{definition}

\begin{definition}
$F:\cN(G)_{\cL} \ra \Set$ is smooth if and only if $F \ra \bullet$ is smooth.
\end{definition}

\begin{lemma}\label{graded} A morphism $f:M \to N$ in $\widehat{\cN}_{\cL}(G)$ is a surjection or an isomorphism if and only if the associated graded morphism $\gr^{\cL} M \to \gr^{\cL} N$ is so. Therefore $f:M \to N$ in $\widehat{\cN}(G)_{\cL}$ is surjective if and only if the induced map \mbox{$t^*_{M/\cL} \to t^*_{N/\cL}$} is surjective, and an endomorphism $f:M \to M$ in  $\widehat{\cN}(G)_{\cL}$ is an automorphism whenever the induced map $t^*_{M/\cL} \to t^*_{M/\cL}$ is the identity.
\end{lemma}

\begin{definition} A map $p:M \to N$ in $\widehat{\cN}(G)_{\cL}$ is essential if for all morphisms $q:M' \to M$, $q$ is surjective whenever $pq$ is.
\end{definition}

From now on, we will assume that $G$ is reductive.

From the above lemma, we deduce:
\begin{lemma}\label{essential}
Let $p:M \to N$ be a surjection in $\widehat{\cN(G)}_{\cL}$. Then
\begin{enumerate}
\item $p$ is essential if and only if the induced map $t^*_{M/\cL} \to t^*_{N/\cL}$ is an isomorphism.
\item If $p$ is a small extension, then $p$ is not essential if and only if $p$ has a section.
\end{enumerate}
\begin{proof}
\begin{enumerate}
\item Assume that $p$ is essential.  Since  $t^*_{M/\cL} \to t^*_{N/\cL}$ is surjective, so must $M \to t^*_{N/\cL}$ be. Using the exactness of $\Rep(G)$ (since $G$ is reductive), this map has a section in $\widehat{\Rep}(G)$. Let its image be $V$, and let $q:L_{\cL}(V) \to M$ be the map determined by the inclusion $V \into M$. Now, $pq$ is surjective, since it induces an isomorphism on cotangent spaces. Therefore $q$ is surjective, so $V \to t^*_{M/\cL}$ is surjective, as required. The converse is immediate.

\item If $p$ is not essential, construct the pro-Lie algebra $L_{\cL}(V)$ as in the previous part, and let $N' \into M$ be its image. Then, by comparing tangent spaces, we see that $N' \to N$ is a surjection. Since $p$ is not essential, $t^*_{M/\cL} \to t^*_{N/\cL}$ is not injective, so $\ker p \to t^*_{M/\cL}$ is non-zero, hence an embedding, since $\ker p$ is absolutely irreducible. The image of $N'$ in $t^*_{M/\cL}$ has zero intersection with  the image of $\ker p$. Therefore $\ker p$ and $N'$ have zero intersection so $N' \by \ker p$ is a sub-pro-Lie algebra of $M$ with the same cotangent space, so $M=N'\by \ker p$, and $N' \cong N$, which gives the section of $p$.
\end{enumerate}
\end{proof}
\end{lemma}    

\begin{definition} A pro-couple $(L,\xi)$ is a hull for $F$ if the induced map $h_L \ra F$ is \'etale.
\end{definition}

\begin{lemma} Suppose $F$ is a functor such that $$F(V\eps\oplus W\eps) \xrightarrow{\sim}F(V\eps)\by F(W\eps)$$ for $V,W \in \Rep(G)$. Then $F(V\eps)$ has a canonical vector space structure, and the tangent space functor 
\begin{eqnarray*}
t_F: \Rep(G) &\to& \Vect\\
V &\mapsto& F(V\eps)
\end{eqnarray*}
is additive.
\end{lemma}

\begin{definition}
Given $F:\cN(G)_{\cL} \ra \Set$, let $N' \ra N$ and $N'' \ra N$ be morphisms in $\cN(G)_{\cL}$, and consider the map:
\begin{enumerate}
\item[$(\dagger)$]  $F(N'\by_N N'') \ra F(N')\by_{F(N)}F(N'').$
\end{enumerate}
We make the following definitions for properties of $F$:
\begin{enumerate}
\item[{\rm (H1)}] $(\dagger)$ is a surjection whenever $N'' \ra N$ is a small extension.
\item[{\rm(H2)}] $(\dagger)$ is a bijection whenever $N=0$ and $N''=V\eps$, for an irreducible $G$-representation $V$.
\item[{\rm(H4)}]$(\dagger)$ is a bijection whenever $N'=N''$ and $N' \ra N$ is a small extension.
\end{enumerate}
\end{definition}

\begin{remark} These conditions are so named for historical reasons, following \cite{Sch}. The missing condition (H3) concerned finite-dimensionality of tangent spaces, which is irrelevant to our (weaker) notion of pro-representability.
\end{remark}

\begin{lemma}\label{lifting} Let $F:\cN(G) \to \Set$ satisfy (H1) and (H2), and $M \to N$ be a semi-small extension in $\widehat{\cN}(G)$. Given a surjection $M \onto M_{\alpha}$, with $M_{\alpha} \in \cN(G)$, let $N_{\alpha}= M_{\alpha}*_M N$. If $\xi \in \hat{F}(N)$ has the property that for all such surjections, the image $\xi_{\alpha} \in F(N_{\alpha})$ of $\xi$ lifts to $F(M_{\alpha})$, then $\xi$ lifts to $\hat{F}(M)$.
\begin{proof}
Since $M \onto N$ is semi-small, $M_{\alpha} \onto N_{\alpha}$ is also. Let \mbox{$I:=\ker(M \to N)$} and \mbox{$I_{\alpha}:=\ker(M_{\alpha} \to N_{\alpha})$.} Observe that we have canonical isomorphisms \mbox{$M_{\alpha} \by I_{\alpha} \cong M_{\alpha}\by_{N_{\alpha}}M_{\alpha}$,} and let $T_{\alpha}$ be the fibre of $F(M_{\alpha} \to F(N_{\alpha}))$ over $\xi_{\alpha}$. Since $F$ satisfies (H2), we have a map
$$
F(M_{\alpha}) \by t_F(I_{\alpha}) \to F(M_{\alpha})\by_{F(N_{\alpha})}F(M_{\alpha}),
$$
so $t_F(I_{\alpha})$ acts on $T_{\alpha}$. From (H1) it follows that this action is transitive (and if $F$ also satisfied (H4) then $T_{\alpha}$ would be a principal homogeneous space under this action). Let $K_{\alpha} \subset t_F(I_{\alpha})$ be the stabiliser of $T_{\alpha}$.

We wish to construct a compatible system $\eta_{\alpha} \in T_{\alpha}$. Let $\eta_{\alpha}' \in T_{\alpha}$ be any element, and assume that we have constructed a compatible system  $\eta_{\beta} \in T_{\beta}$, for all strict epimorphisms $M_{\alpha} \onto M_{\beta}$. Then $\eta_{\beta} = v_{\beta}(\eta_{\alpha}')$, for a unique $v_{\beta} \in  t_F(I_{\beta})/K_{\beta}$. Observe that
$$
t_F(I_{\alpha})/K_{\alpha} \to \lim_{\substack{\longleftarrow \\ \beta}} t_F(I_{\beta})/K_{\beta}
$$
is surjective, and that the $v_{\beta}$ form an element of the right-hand side. Lift to \mbox{$v \in t_F(I_{\alpha})/K_{\alpha}$,} and let $\eta_{\alpha}:=  v(\eta_{\alpha}')$. The construction proceeds inductively (since every poset can be enriched to form a totally ordered set, and we have satisfied the hypotheses for transfinite induction).
\end{proof}
\end{lemma}

\begin{proposition} Let $(L, \xi)$, $(L',\xi')$ be hulls of $F$. Then there exists an isomorphism $u:L \ra L'$ such that $F(u)(\xi)=\xi'$.
\begin{proof} We wish to lift $\xi$ and $\xi'$ to $u\in h_{L'}(L)$ and $u' \in h_L(L')$. We may apply Lemma \ref{lifting} to the the functor $h_{L'}$ and the  successive semi-small extensions $L_{n+1} \to L_n$, where $L_n=L/(\Gamma_n L)$. By smoothness of $h_{L'} \to F$, we obtain successive lifts of $\xi_n \in F(L_n)$ to $u_n \in h_{L'}(L_n)$. Let $u = \Lim u_n$, and construct $u'$ similarly. We therefore obtain $u:(L',\xi) \to (L,\xi)$ and $u': (L,\xi)\to(L',\xi)$, inducing identity on cotangent spaces. Therefore $uu'$ induces the identity on $t^*_{L/\cL}$, so is an automorphism, by Lemma \ref{graded}.
\end{proof}
\end{proposition}

\begin{proposition}
\begin{enumerate}
\item Let $M \ra N$ be a morphism in $\widehat{\cN(G)}_{\cL}$. Then $h_N \ra h_M$ is smooth if and only if $N\cong L_M(V)$, for some $V \in \Rep(G)$.
\item If $F \ra E$ and $E \ra H$ are smooth morphisms of functors, then the composition $F \ra H$ is smooth.
\item If $u:F \ra E$ and $v:E \ra H$ are morphisms of functors such that u is surjective and $vu$ is smooth, then $v$ is smooth.
\end{enumerate}
\begin{proof} If $N=M*L(V)$, then 
\mbox{$
h_N(\g)= h_M(\g)\by \Hom(V,\g),
$}
 which is smooth over $h_M$, since $\Rep(G)$ is an exact category. Conversely, assume that $h_N \ra h_M$ is smooth, and let $V=t^*N/M$. Let $N'=L_M(V)$, and observe that, by choosing a lift of $V$ to $N$, we obtain a morphism $f:N' \to N$, inducing an isomorphism on relative cotangent spaces.

Let $N'_n= N'/(\Gamma_{n+1} N')$, and observe that $N'_n \to N'_{n-1}$ is a semi-small extension, as is $N'_1 \to t^*_{N'/M}$. We have a canonical map $N \to  t^*_{N'/M}$ arising from the isomorphism $t^*_{N'/M} \cong  t^*_{N/M}$. Since $h_N$ satisfies (H1) and (H2), we may now apply Lemma \ref{lifting} to construct a map $g:N \to N'$ lifting this. Therefore the compositions $fg$ and $gf$ induce the identity on cotangent spaces, so are isomorphisms, and $N'\cong N$, as required.

The remaining statements follow by formal arguments. 
\end{proof}
\end{proposition}

\begin{remark}
For the proposition above, it is essential that $G$ be reductive, since we need the exactness of $\Hom$ on $\Rep(G)$ to ensure that $L(V)$ is smooth.
\end{remark}

\begin{theorem}\label{nSch}
\begin{enumerate}
\item $F$ has a hull if and only if $F$ has properties (H1) and (H2).
\item $F$ is moreover pro-representable if and only if $F$ has the additional property (H4).
\end{enumerate}
\begin{proof}
This is essentially \cite{Sch}, Theorem 2.11, \emph{mutatis mutandis}.
Since $t_F$ defines a left-exact functor on $\Rep(G)$, by (H2), let it be pro-represented by  $W \in \widehat{\Rep}(G)$, and let $\fh=L_{\cL}(W)$. The hull $\g$ will be a quotient of $\fh$, which will be constructed as an inverse limit of semi-small extensions of $W$. Let $\g_2=W$, and $\xi_2 \in \hat{F}(W)$ the canonical element corresponding to the pro-representation of $t_F$. Assume that we have constructed $(\g_q, \xi_q)$. We wish to find a semi-small extension 
$\g_{q+1} \to \g_q$, maximal among those quotients of $\fh$ which admit a lift of $\xi_q$.

Given quotients $M,N \in \cN(G)$ of $\fh$, we write 
$$
M\wedge N:=M*_{\fh}N,\quad \text{ and } \quad M\vee N:= M\by_{M\wedge N} N.
$$
 Note that $M\wedge N$ is then maximal among those Lie algebras which are dominated by both $M$ and $N$, while $M\vee N$ is minimal among those quotients of $\fh$ which dominate both $M$ and $N$. Next, observe that a set $\{M_{\alpha}\in \cN(G)\}_{\alpha \in I}$ of quotients of $\fh$ corresponds to a quotient $\fh \onto \Lim_{\alpha} M_{\alpha}$ of $\fh$ if and only if the following two conditions hold:
\begin{enumerate}
\item[{\rm (Q1)}] If $M_{\alpha} \onto N$ is a surjection, for any $\alpha \in I$ and any $N \in \cN(G)$, then $N \cong M_{\beta}$, for some $\beta \in I$.
\item[{\rm (Q2)}] Given any $\alpha, \beta \in I$, $M_{\alpha} \vee M_{\beta} =M_{\gamma}$, for some $\gamma \in I$. 
\end{enumerate}

We will now form such a set of quotients by considering those $\fh \onto M$ satisfying:
\begin{enumerate}
\item $M \onto M*_{\fh}\g_q$ is a semi-small extension.
\item The image of $\xi_q$ in $F(M*_{\fh}\g_q)$ lifts to $F(M)$.
\end{enumerate}
It is immediate that this set satisfies (Q1). To see that it satisfies (Q2), take quotients $M,N$ satisfying these conditions. It is clear that $M\vee N \onto (M\vee N)*_{\fh} \g_q$ is a semi-small extension, since $M\vee N$ is a sub-Lie algebra of $M\by N$. To see that $\xi_q$ lifts to $F(M\vee N)$, let $x \in F(M), y \in F(N)$ be lifts of $\xi_q$. Now, as in the proof of Lemma \ref{lifting}, the fibre of $F(N)$ over $\xi_q$ surjects onto the fibre of $F(M\wedge N)$ over $\xi_q$. Therefore, we may assume that $x$ and $y$ have the same image in $F(M\wedge N)$. Now (H1) provides the required lift:
$$
F(M\vee N) \onto F(M) \by_{F(M\wedge N)}F(N).
$$

Let $\g_{q+1} \in \widehat{\cN}(G)$ be the quotient defined by this collection. By Lemma \ref{lifting}, it follows that $\xi_q$ lifts to $\hat{F}(\g_{q+1})$. Let $\g:=\Lim \g_q$, with $\xi:=\Lim \xi_q$. 

It remains to show that this is indeed a hull for $F$. By construction, $h_{\g} \to F$ is unramified, so we must show it is smooth. Let $p:(N',\eta') \to (N,\eta)$ be a morphism of couples in $\cN(G)$, with $p$ a small extension, $N=N'/I$, and assume we are given $u:(\g,\xi) \to (N,\eta)$. We must lift $u$ to a morphism $(\g,\xi) \to (N'\eta')$. It will suffice to find a morphism $u':\g \to A'$ such that $pu'=u$, since we may then use transitivity of the action of $t_F(I)$ on $F(p)^{-1}(\eta)$.

For some $q$, $u$ factors as $(\g,\xi) \to (\g_q,\xi_q) \to (N,\eta)$. By smoothness of $\fh$, we may choose a morphism $w$ making the following diagram commute. We wish to construct the morphism $v$:
$$
\xymatrix{ \fh \ar[r]^-{w} \ar[d] & \g_q\by_NN' \ar[d]^{\mathrm{pr}_1}\\
\g_{q+1} \ar@{--}[ur]^-{v} \ar[r]& \g_q.
} 
$$

If the small extension $\mathrm{pr}_1$ has a section, then $v$ obviously exists. Otherwise, by Lemma \ref{essential}, $\mathrm{pr}_1$ is essential, so $w$ is a surjection. (H1) then provides a lift of $\xi_q$ to $F(\g_q\by_NN')$, so by the construction of $\g_{q+1}$, $w$ factors through $\g_{q+1}$, and so $v$ must exist. This completes the proof that $h_{\g} \to F$ is a hull. 

If $F$ also satisfies (H4), then we may use induction on the length of $N$ to show that $h_{\g}(N) \cong F(N)$, using the observation in the proof of Lemma \ref{lifting} that all non-empty fibres over small extensions $I \to N' \to N$ are principal homogeneous $t_F(I)$-spaces. 

Necessity of the conditions follows by a formal argument.
\end{proof}
\end{theorem}

\begin{definition}
$F:\cN(G)_{\cL} \ra \Set$ is homogeneous if 
$$
\eta : F(N'\by_N N'')\ra F(N')\by_{F(N)}F(N'')
$$
is an isomorphism for every $N' \twoheadrightarrow N$.

Note that a homogeneous functor satisfies conditions (H1), (H2) and (H4).
\end{definition}

\begin{definition}
$F:\cN(G)_{\cL} \ra \Set$ is a deformation functor if:
\begin{enumerate}
\item $\eta$ is surjective whenever $N' \twoheadrightarrow N$.
\item $\eta$ is an isomorphism whenever $N=0$.
\end{enumerate}
Note that a deformation functor satisfies conditions (H1) and (H2).
\end{definition}

The following results are adapted from\cite{Man}:

\begin{definition} Given $F:\cN_{\cL,k} \to \Set$, an obstruction theory $(O,o_e)$ for $F$ consists of
an  additive functor $O:\Rep(G) \to \Vect$, the obstruction space, together with obstruction maps $o_e: F(N) \to O(I)$ 
 for each small extension
$$
e: 0 \to I \to L \to N \to 0,
$$
such that:
\begin{enumerate}
\item If $\xi \in F(N)$ can be lifted to $F(L)$ then $o_e(\xi)=0$.
\item For every morphism $\alpha:e \to e'$ of small extensions, we have \mbox{$o_{e'}(\alpha(\xi))= O(\alpha)(v_e(\xi))$,} for all $\xi \in F(N)$.
\end{enumerate}

An obstruction theory $(O,o_e)$ is called complete if  $\xi \in F(N)$ can be lifted to $F(L)$ whenever $o_e(\xi)=0$.
\end{definition}

\begin{proposition}\label{nSSC} (Standard Smoothness Criterion) Given $\phi :F \ra E$, with \mbox{$(O,o_e) \xrightarrow{\phi'}(P,p_e)$} a compatible morphism of obstruction theories, if $(O,o_e)$ is complete, $O \xrightarrow{\phi'} P$ injective, and $t_F \ra t_E$ surjective, then $\phi$ is smooth.
\begin{proof} \cite{Man}, Proposition 2.17.
\end{proof}
\end{proposition}

For functors $F:\cN(G)_{\cL} \ra \Set$ and $E:\cN(G)_{\cL} \ra \Grp$, we say that $E$ acts on $F$ if we have a functorial group action $E(N) \by F(N) \xra{*} F(N)$, for each $N$ in $\cN(G)$.  The quotient functor $F/E$ is defined by $(F/E)(N)=F(N)/E(N)$.

\begin{proposition}\label{nMan1} If $F:\cN(G)_{\cL} \ra \Set$, a deformation functor, and $E:\cN(G)_{\cL} \ra \Grp$ a smooth deformation functor, with $E$ acting on $F$, then $D:=F/E$ is a
deformation functor, and if $\nu :t_E \ra t_F$ denotes $h \mapsto h*0$, then $t_D=\coker\nu$, and the universal obstruction theories of $D$ and $F$ are isomorphic.
\begin{proof} \cite{Man}, Lemma 2.20.
\end{proof}
\end{proposition}

\begin{proposition}\label{nMan2} For $F:\cN(G)_{\cL} \ra \Set$  homogeneous, and $E:\cN(G)_{\cL} \ra \Grp$ a deformation functor, given $a,b \in F(L)$, define $\Iso(a,b): \cN_{L,k} \ra \Set$  by 
$$ 
\Iso(a,b)(L \xrightarrow{f}N)=\{g \in E(N) | g*f(a)=f(b)\}.
$$ 
Then $\Iso(a,b)$ is a deformation functor, with tangent space $\ker\nu$ and, if $E$ is moreover smooth, complete obstruction space $\coker\nu=t_D$.
\begin{proof} \cite{Man}, Proposition 2.21.
\end{proof}
\end{proposition}

\begin{proposition}\label{nMan3} If $E,E'$ are smooth deformation functors, acting on $F,F'$ respectively, with $F,F'$ homogeneous, $\ker \nu \ra \ker \nu'$ surjective, and $\coker\nu\ra\coker\nu'$ injective, then $F/E \ra F'/E'$ is injective.
\begin{proof} \cite{Man}, Corollary 2.22.
\end{proof}
\end{proposition}

This final result does have an analogue in \cite{Man}, but proves extremely useful:

\begin{corollary}\label{nkeyhgs}
If $F:\cN(G)_{\cL} \ra \Set$ and $E:\cN(G)_{\cL} \ra \Grp$ are deformation functors, with $E$ acting on $F$, let $D:=F/E$, then:
\begin{enumerate}
\item If $E$ is smooth, then $\eta_D$ is surjective for every $M \onto N$ (i.e. $D$ is a deformation functor).
\item If $F$ is homogeneous and $\ker \nu=0$, then $\eta_D$ is injective for every $M \onto N$.
\end{enumerate}
Thus, in particular, $F/E$ will be homogeneous if $F$ is homogeneous, $E$ is a smooth deformation functor and $\ker \nu=0$.
%\begin{proof}
%As for \cite{sdc} Theorem \ref{sdc-keyhgs}.
%\end{proof}
\end{corollary}

To summarise the results concerning the pro-representability of the quotient \mbox{$D=F/E$,} we have:
\begin{enumerate}
\item If $F$ is a deformation functor and $G$ a smooth deformation functor,  then $D$ has a hull.
\item If $F$ is homogeneous and $E$ a smooth deformation functor, with $\ker \nu =0$, then $D$ is pro-representable.
\end{enumerate}

\section{Twisted differential graded algebras}\label{tdga}

Throughout this section, we will adopt the conventions of \cite{tannaka} concerning tensor categories. In particular, the associativity isomorphisms will be denoted
$$
\phi_{UVW}: U\ten (V\ten W) \to (U\ten V)\ten W,
$$
and the commutativity isomorphisms
$$
\psi_{UV}: U \ten V \to V \ten U.
$$

\begin{definition} The category $\DGVect$ of graded real vector spaces $\bigoplus_{i \ge 0} V^i$ is a tensor category, with the obvious tensor product 
$$
(U\ten V)^n = \bigoplus_{i+j=n} U^i \ten V^j,
$$
with differential $d|_{U^i \ten V^j}=  d_{U^i}\ten \id + (-1)^j \id \ten d_{V^j}$. The associativity map is the obvious one, while the commutativity map is
\begin{eqnarray*}
\psi_{UV}:U \ten V &\to& V \ten U,\\
u\ten v &\mapsto& (-1)^{ij} v\ten u,
\end{eqnarray*}
for $u \in U^i, v \in V^j$.
\end{definition}

\begin{definition} A (real) DGA over a tensor category $\C$ is defined to be an additive functor $A: \C  \to \DGVect$, equipped with a multiplication
$$
\mu_{UV}: A(U) \ten A(V) \to  A(U\ten V),
$$
functorial in $U$ and $V$, such that
\begin{enumerate} 
\item Associativity. The following diagram commutes:
$$
\xymatrix{
AU \ten (AV \ten AW) \ar[d]_{\phi} \ar[r]^{\id\ten \mu} & AU\ten A(V\ten W) \ar[r]^{\mu} &A(U\ten(V\ten W)) \ar[d]^{A\phi}\\
(AU\ten AV)\ten AW \ar[r]^{\mu \ten \id}  &A(U\ten V)\ten AW \ar[r]^{\mu} &A((U\ten V)\ten W).
}
$$

\item  Commutativity. The following diagram commutes:
$$
\xymatrix{
AU\ten AV \ar[r]^{\mu} \ar[d]_{\psi} &A(U\ten V) \ar[d]^{A\psi} \\
AV\ten AU \ar[r]^{\mu} & A(V\ten U).
}
$$
\end{enumerate}
\end{definition}

\begin{remark}
Note that, if we take $\C$ to be the category of finite-dimensional complex vector spaces, then  giving  a DGA $A$ over $\C$ is equivalent to giving the differential graded algebra $A(\R)$, which motivates the terminology.
\end{remark}

\begin{definition} We say that a DGA $A$ over $\C$ is flat if for every exact sequence \mbox{$0 \to U \to V \to W \to 0$} in $\C$, the sequence $0 \to AU \to AV \to AW \to 0$ is exact.
\end{definition}

\begin{lemma} Given a finite-dimensional Lie algebra $L$ with a $G$-action, the graded vector space $A(L)$ has the natural structure of a differential graded Lie algebra.
\begin{proof} It suffices to define the Lie bracket. Let it be the composition 
$$
A(L) \ten A(L) \xra{\mu} A(L\ten L) \xra{A([,])} A(L).
$$
the associativity, commutativity  and compatibility axioms are easily verified.
\end{proof}
\end{lemma}

\begin{definition} Given a flat DGA $A$ over $\Rep(G)$,   the Maurer-Cartan functor 
\mbox{$\mc_A:\cN(G) \ra \Set$} is defined by 
$$ \mc_A(N)=\{x \in A(N)^1 |dx+\half[x,x]=0\}.$$
\end{definition}

Observe that for $\omega \in A(N)^1$, 
$$
d\omega +\half[\omega,\omega]=0 \Rightarrow (d+\ad_{\omega})\circ(d+\ad_{\omega})=0,
$$ 
so $(A(N), [,], d+\ad_{\omega})$ is a DGLA.

\begin{definition}
Define the gauge functor $\ga:\cN(G) \ra \Grp$ by 
$$
\ga(N)=\exp(A(N)^0),
$$ 
noting that nilpotence of $N$ implies nilpotence of $A(N)^0$.
\end{definition}

We may now define the  DGLA $(A(N))_d$ as in \cite{Man}:
$$
(A(N))_d^i=\left\{\begin{matrix} (A(N))^1\oplus \R d	& i=1\\
                           (A(N))^i 			& i \ne 1, \end{matrix}\right.
$$
 with 
$$
d_d(d)=0,\quad [d,d]=0,\quad [d,a]_d=da,\quad \forall a \in (A(N)).
$$ 

\begin{lemma} 
 $\exp(A(N)^0)$ commutes with $[,]$ when acting  on $(A(N))_d$ via the adjoint action.
\end{lemma}

\begin{corollary}
Since $\exp(A(N)^0)$ preserves $(A(N)^1)+d \subset (A(N))_d$ under the adjoint action, and 
$$
x \in \mc_A(N) \iff [x+d,x+d]=0,
$$
 the adjoint action of $\exp(A(N)^0)$ on $A(N)^1+d$ induces an action of $\ga(N)$ on $\mc_A(N)$, which we will call the gauge action.
\end{corollary}

\begin{definition} $\defa=\mca/\ga$, the quotient being given by the gauge action \mbox{$\alpha(x)= \ad_{\alpha}(x+d)-d$.}  Observe that $\ga$  and $\mca$ are homogeneous.  Define the deformation groupoid $\mathfrak{Def}_A$ to have objects $\mca$, and morphisms  given by $\ga$.
\end{definition}

Now, 
$
t_{\ga}(V)=A^0(V) ,
$ 
and 
$
t_{\mca}(V)=\z^1(A(V) ),
$
 with action 
\begin{eqnarray*}
t_{\ga}\by t_{\mca} &\ra& t_{\mca};\\ 
(b,x) &\mapsto& x+db, \text{ so }
\end{eqnarray*}
$$
 t_{\defa}(V)=\H^1(A(V)). 
$$  

\begin{lemma} $\H^2(A)$ is a complete obstruction space for
$\mca$.
\begin{proof} 
Given a small extension 
$$
e:0 \ra I \ra N \ra M\ra 0,
$$ 
and $x \in \mca(M)$, lift $x$ to $\tilde{x} \in A^1(N)$, and let 
$$
h=d\tilde{x} +\half [\tilde{x},\tilde{x}] \in A^2(N).
$$
 In fact, $h \in A^2(I)$, as $dx+\half[x,x]=0$. 

Now,
$$
dh=d^2\tilde{x}+[d\tilde{x},\tilde{x}]= [h-\half [\tilde{x},\tilde{x}],\tilde{x}]=[h,\tilde{x}]=0,
$$ 
since $[[\tilde{x},\tilde{x}],\tilde{x}]=0$ and $[I,N]=0$. Let 
$$
o_e(x)=[h]\in \H^2(A(I)).
$$
This is well-defined: if $y=\tilde{x}+z$, for $z \in A^1(K)$, then
$$
dy +\half [y,y]=d\tilde{x}+dz+\half[\tilde{x},\tilde{x}]+\half[z,z]+[\tilde{x},z]=h+dz,
$$ 
as $[I,N]=0$.

This construction is clearly functorial, so it follows that $(\H^2(A),o_e)$ is a complete obstruction theory
for $\mca$.
\end{proof}
\end{lemma}

Now Proposition \ref{nMan1} implies the following:

\begin{theorem} $\defa$ is a deformation functor, $t_{\defa} \cong \H^1(A)$, and $\H^2(A)$ is a complete obstruction theory for $\defa$.
\end{theorem}

The other propositions of Section \ref{nilplie} can be used to prove:

\begin{theorem}\label{nqis}
If $\phi :A \ra B$ is a morphism of DGAs over $\Rep(G)$, and 
$$
\H^i(\phi):\H^i(A) \ra \H^i(B)
$$ 
are the induced maps on cohomology, then:
\begin{enumerate}
\item If $\H^1(\phi)$ is bijective, and $\H^2(\phi)$ injective, then $\defa \ra \ddef_B$ is \'etale.
\item If also $\H^0(\phi)$ is surjective, then  $\defa \ra \ddef_B$ is an isomorphism.
\item Provided condition 1 holds, $\mathfrak{Def}_A \to \mathfrak{Def}_B$ is an equivalence of functors of  groupoids if and only if $\H^0(\phi)$ is an isomorphism.
\end{enumerate}
\begin{proof}
\cite{Man}, Theorem 3.1, mutatis mutandis.
\end{proof}
\end{theorem}

\begin{theorem}
If $\H^0(A)=0$, then $\defa$ is homogeneous.
\begin{proof}
Proposition \ref{nkeyhgs}.
\end{proof}
\end{theorem}

Thus, in particular, a quasi-isomorphism of DGAs gives an isomorphism of deformation functors and of deformation groupoids.

\begin{remark}\label{toen} The category of DGAs over $\Rep(G)$ is, in fact, equivalent to the category of $G$-equivariant differential graded algebras. Given a DGA $A$ over $\Rep(G)$, we consider the structure sheaf $O(G)$ of $G$, regarded as a $G$-representation via the left action. Then $O(G) \in \ind(\Rep(G))$, and  we therefore set $B=A(O(G))$, which has a DGA structure arising from the algebra structure on $O(G)$, and a $G$-action given by the right action on $O(G)$. Conversely, given a $G$-equivariant DGA $B$, we define $A(V):=B\ten^G V $, the subspace of $G$-invariants of $B\ten V$. 

By Lemma \ref{dual}, the vector space $O(G)\ten^G V$  is isomorphic to $V$, with the $G$-action on $O(G)$ coming from the right action of $G$. This implies that the functors above define an equivalence. This equivalence will mean that the twisted DGA considered in Section \ref{hodge} is a model for the schematic homotopy type considered in \cite{KTP}.
\end{remark}

\section{Relative Malcev completions}\label{proalg}

\begin{definition} Given  a  group $\Gamma$ with  a representation $\rho_0:\Gamma \to G$ to a reductive real  pro-algebraic group, define the functor
$$
\fR_{\rho_0}:\cN(G) \to \Grpd
$$
of deformations of $\rho_0$ so that the objects of $\fR_{\rho_0}(\fu)$ are representations
$$
\rho: \Gamma \to \exp(\fu) \rtimes G
$$
lifting $\rho_0 $, and isomorphisms are given by the conjugation action of the unipotent group $\exp(\fu)$ on $\exp(\fu) \rtimes G$. Explicitly, $u \in \exp(\fu)$ maps $\rho$ to $ u\rho u^{-1}$.
\end{definition}

\begin{lemma}
\begin{enumerate}
\item The functor $R_{\rho_0}$ of objects of $\fR_{\rho_0}$ is a deformation functor, with tangent space $V \mapsto \H^1(\Gamma,\rho_0^{\sharp}V)$ and obstruction space  $V \mapsto \H^2(\Gamma,\rho_0^{\sharp}V)$. 

\item Given $\omega , \omega' \in \fR_{\rho_0}(\g)$, the functor on $\cN(G)_{\g}$ given by
$$
\fu \mapsto \Iso_{\fR_{\rho_0}(\fu)}(\omega, \omega')
$$
is homogeneous, with tangent space
$$
V \mapsto \H^0(\Gamma,\rho_0^{\sharp}V)
$$
and obstruction space
$$
V \mapsto \H^1(\Gamma,\rho_0^{\sharp}V).
$$
\end{enumerate}
\end{lemma}

\begin{proposition}\label{gammahat} Let $\Gamma^{\alg}$ be the pro-algebraic completion of $\Gamma$, and $\rho_0:\Gamma \to \Gamma^{\red}$ its reductive quotient. Then the Lie algebra $ \cL(\Ru(\Gamma^{\alg}))$  of the pro-unipotent radical of $\Gamma^{\alg}$, equipped with its $\Gamma^{\red}$-action as in Theorem \ref{radical}, is a hull for the functor $R_{\rho_0}$.
\begin{proof}
By definition,
\begin{eqnarray*}
R_{\rho_0}(U)&=&\Hom(\Gamma, U \rtimes \Gamma^{\red})_{\rho_0}/U\\
&=& \Hom(\Gamma^{\alg}, U \rtimes \Gamma^{\red})_{\rho_0}/U.
\end{eqnarray*}
If we now fix a Levi decomposition $\Gamma^{\alg}\cong \Ru(\Gamma^{\alg})\rtimes \Gamma^{\red}$, we may rewrite this as
$$
 \Hom(\Ru(\Gamma^{\alg})\rtimes \Gamma^{\red}, U \rtimes \Gamma^{\red})_{\rho_0}/U.
$$

There is  a natural map
$$
f:\Hom_{\hat{\cN}(\Gamma^{\red})}(\Ru(\Gamma^{\alg}), U) \to \Hom(\Ru(\Gamma^{\alg})\rtimes \Gamma^{\red}, U \rtimes \Gamma^{\red})_{\rho_0}/U,
$$
and we need to show that this map is surjective, and an isomorphism on tangent spaces. 

For surjectivity, take 
$$
\rho:\Ru(\Gamma^{\alg})\rtimes \Gamma^{\red}\to U \rtimes \Gamma^{\red}
$$
lifting $\rho_0$. Since $\Gamma^{\red}$ is reductive, $\rho(\Gamma^{\red})\le U \rtimes \Gamma^{\red}$ must be reductive. But the composition
$$
\rho(\Gamma^{\red}) \into U \rtimes \Gamma^{\red} \onto \Gamma^{\red}
$$
is a surjection, and $\Gamma^{\red}$ is also the reductive quotient of $U \rtimes \Gamma^{\red}$, so $\rho(\Gamma^{\red})$ is a maximal reductive subgroup. By the Levi decomposition theorem, maximal reductive subgroups are conjugate under the action of $U$, so     there exists $u \in U$ such that $\ad_u\rho(\Gamma^{\red})=\Gamma^{\red}$. Now we may replace $\rho$ by $\ad_u\rho$, since they  define the same element of $R_{\rho_0}(U)$. As $\ad_u\rho$ preserves $\Gamma^{\red}$, its restriction  $\Ru(\Gamma^{\alg})\to U$ is $\Gamma^{\red}$-equivariant, so $\ad_u\rho$ lies in the image of $f$.

To see that $f$ induces an isomorphism on tangent spaces, we need to show that it 
is injective whenever $U$ is abelian. This is immediate, since the conjugation action of $U$ on $U$ is then trivial.
\end{proof}
\end{proposition}

\begin{remark} In the terminology of \cite{malcev}, $\Gamma^{\alg} \to \Gamma^{\red}$ is the relative Malcev completion of the representation $\Gamma \to \Gamma^{\red}$, so we can regard this section as studying Malcev completions of arbitrary Zariski-dense reductive representations.
\end{remark}

\begin{definition}\label{sharp} Given a homomorphism $\theta:G \to H$ of algebraic  groups, with $H$ reductive, define \mbox{$\theta_{\sharp}:\widehat{\cN}(G) \to \widehat{\cN}(H)$} to be left adjoint to the restriction map \mbox{$\theta^{\sharp}: \widehat{\cN}(H)\to \widehat{\cN}(G)$,} so that
$$
\Hom_{\widehat{\cN}(G)}(\theta_{\sharp}L,N) \cong \Hom_{\widehat{\cN}(H)}(L, \theta^{\sharp} N).
$$
This left adjoint must exist, since the functor on the right satisfies Schlessinger's conditions.  
\end{definition}

\begin{lemma}\label{semidirect} If $\Gamma = \Delta \rtimes \Lambda$,  such that the adjoint action of $\Lambda$ on the pro-unipotent completion $\Delta \ten \R$ is reductive, then
$$
\bar{\rho}_{\sharp}\Ru(\Gamma^{\alg}) \cong (\Delta \ten \R)\by \Ru(\Lambda^{\alg}) \in \exp(\widehat{\cN}(\Lambda^{\red})),
$$
where we write $\rho$ for the  composition $\Gamma \to \Lambda^{\red}$, and $\bar{\rho}$ for the  quotient  representation $\Gamma^{\red} \to \Lambda^{\red}$. 
\begin{proof}
We use the fact that $\bar{\rho}_{\sharp}\Ru(\Gamma^{\alg})$ pro-represents the functor \mbox{$U \mapsto \Hom(\Gamma,  U \rtimes \Lambda^{\red})_{\rho}$,} for $U \in \exp(\cN(\Lambda^{\red}))$. 

A homomorphism $\Gamma \to U \rtimes \Lambda^{\red}$ lifting $\rho$ gives rise to a map $ \Delta\ten \R \to U$, since $U$ is unipotent. 
It is then clear that
$$
\bar{\rho}_{\sharp}\Ru(\Gamma^{\alg}) \cong ( \Delta\ten \R) \rtimes \Ru(\Lambda^{\alg}).
$$
Finally, observe that $\Ru(\Lambda^{\alg})$ acts trivially on $\Delta\ten \R$, since the action is reductive.
\end{proof}
\end{lemma}

\section{Principal homogeneous spaces}\label{torsor}

Fix a connected differentiable manifold $X$. Let $\redpi $ be the reductive quotient of the pro-algebraic real completion of $\pi_1(X,x)$, so that $\Rep(\redpi)$ can be regarded as the category of real semisimple $\pi_1(X,x)$-representations. Given such a representation $V$, let $\vv$ denote the corresponding semisimple local system.

\begin{definition}
We may then define a DGA over $\Rep(\redpi)$ by 
$$
A(V):=\Gamma(X,\vv \ten\sA^{\bullet}),
$$
 where $\sA^{\bullet}$ is the sheaf of real $\CC^{\infty}$ forms on $X$. The multiplication is given by 
\begin{eqnarray*}
A(V)\ten A(W) &\cong& \Gamma(X\by X, p_1^*(\vv \ten\sA^{\bullet})\ten p_2^*(\ww \ten\sA^{\bullet}))\\
&\xra{\Delta^{\sharp}}& \Gamma(X,(\vv \ten\sA^{\bullet})\ten (\ww \ten \sA^{\bullet}))\\
&\to& \Gamma(X, (\vv \ten \ww)\ten  \sA^{\bullet}),
\end{eqnarray*}
the first map being the K\"unneth isomorphism, where $\Delta:X \to X\by X$ is the diagonal map, and $p_1,p_2: X\by X \to X$ the projection maps. The final isomorphism is the composition of the multiplication on $\sA^{\bullet}$ with the relevant associativity and commutativity isomorphisms.
\end{definition}

The aim of this section is to prove that the groupoids $\Def_A(\g)$ are functorially equivalent to the groupoid of $\exp(\mathpzc{g})$-torsors, where $\mathpzc{g}$ is the sheaf of Lie algebras associated to $\g$. 

\begin{definition} Given a locally constant  sheaf $\bG$ of groups on $X$, define $\fB(\bG)$, the category of $\bG$-torsors (or principal homogeneous $\bG$-spaces) to consist of sheaves of sets $\Bu$ on $X$, together with a multiplication $\bG \by \bB \to \bB$ such that  $g\cdot(h\cdot b)=(gh)\cdot b$, and the stalks $\bB_x$ are isomorphic (as $\bG_x$-spaces) to $\bG_x$.
\end{definition}

\begin{lemma} There is a canonical morphism $\Bu:\Def_A(\g) \to \fB(\exp(\mathpzc{g}))$, functorial in $\g \in \cN(\redpi)$.
\begin{proof}
Given $\omega \in \mca$, let
$$
\Bu_{\omega}:=D^{-1}(\omega),
$$
where
\begin{eqnarray*}
D: \exp(\mathpzc{g} \ten \sA^0) &\to& \mathpzc{g} \ten \sA^1\\
\alpha &\mapsto& d\alpha\cdot\alpha^{-1}.
\end{eqnarray*}
Then $\Bu_{\omega}$ is a principal $\exp(\g)$-sheaf on $X$.
\end{proof}
\end{lemma}

\begin{lemma} \begin{enumerate} 
\item The functor  $\g \mapsto B(\exp(\mathpzc{g}))$, the set of isomorphism classes of $\fB(\exp(\mathpzc{g}))$, is a deformation functor with tangent space 
$$
V \mapsto \H^1(X,\vv),
$$
and obstruction space
$$
V \mapsto \H^2(X, \vv).
$$

\item Given $\omega , \omega' \in \fB(\exp(\mathpzc{g}))$, the functor on $\cN(\redpi)_{\g}$ given by
$$
\fh \mapsto \Iso_{\fB(\exp(\mathpzc{h}))}(\omega, \omega')
$$
is homogeneous, with tangent space
$$
V \mapsto \H^0(X,\vv)
$$
and obstruction space
$$
V \mapsto \H^1(X, \vv).
$$
\end{enumerate}
\begin{proof} 
Take an cover $\{U_i\}$ of $X$ by  open discs. Then a $\vv$-torsor $\Bu$ is determined by fixing isomorphisms $\vv|_{U_i}\cong \Bu|_{U_i}$ and specifying transition maps in $\vv_{U_i\cap U_j}$ satisfying the cocycle condition. The result follows by considering isomorphism classes of these data.
\end{proof}
\end{lemma}

\begin{theorem}\label{vbeq} The functor $\Bu$ is an equivalence of groupoids.
\begin{proof}
We begin by proving essential surjectivity. The morphism $\ddef_A(\g) \to B(\exp(\mathpzc{g}))$ induces an isomorphism on tangent and obstruction spaces $\H^i(X,\vv)$, so is \'etale (by Proposition \ref{nSSC}).
Now, $\Iso(\omega,\omega') \to \Iso(\Bu_{\omega},\Bu_{\omega'})$ on $\cN(G)_{\g}$ is similarly \'etale, so must be an isomorphism, both functors being pro-representable.
\end{proof}
\end{theorem}

We will look at an algebraic interpretation of the groupoids we have been considering.

\begin{lemma}\label{repeq} If $X$ is a connected differentiable manifold and $\Gamma =\pi_1(X,x)$ is its fundamental group, then  there is a canonical equivalence of groupoids 
$$
\Bu: \fR_{\rho_0}(\g) \to \fB(\exp(\mathpzc{g})),
$$
for  $\rho_0:\pi_1(X,x) \to \redpi$, and $\g \in \cN(\redpi)$.
\begin{proof} Let $\tilde{X} \xra{\pi} X$ be the universal covering space  of X, on which $\Gamma$ acts. Then, associated to any representation $\rho: \Gamma \to H$, we have the $H$-torsor 
$$
\Bu_{\rho,H}:= (\pi_*H)^{\Gamma,\rho}.
$$
Associated to any  $\rho:\Gamma \to \exp(\g) \rtimes G$ lifting $\rho_0$, we have a representation $\rho: \Gamma \to H$, where $H=\exp(\g) \rtimes \Gamma$. This gives rise to the $H$-torsor $\Bu_{\rho,H}$.   Let $\Bu_{\rho}:=\Bu_{\rho,H}/\Gamma$ be the quotient sheaf under the $\Gamma$-action (using $\Gamma \le H$). It follows that this is an $\exp(\mathpzc{g})$-torsor.

Finally to see that $\Bu$ defines an equivalence, observe that the maps on tangent and obstruction spaces are
$$
\H^i(\Gamma,V) \to \H^i(X,\vv),
$$
which are isomorphisms for $i=0,1$, and injective for $i=2$. The equivalence then follows from Proposition \ref{nSSC}. 
\end{proof}
\end{lemma}

\section{Hodge theory}\label{hodge}

Let $X$ be a compact connected K\"ahler manifold, and  let $\redpi$ be the reductive pro-algebraic completion of $\pi_1(X,x)$. Recall that the DGA $A$ is defined over $\Rep(\redpi)$ by
$$
A(V):=\Gamma(X, \vv\ten \sA^{\bullet}).
$$

Since $\Rep(\redpi)$ is an exact category, note that all DGAs over $\Rep(\redpi)$ are flat.

\begin{theorem}\label{formal}
The DGA $A$ is formal, i.e. weakly equivalent to its cohomology DGA.
\begin{proof}
We have an operator $\dc=J^{-1}dJ$ on $A$, where $J$ is the complex structure. This satisfies $d\dc+\dc d=0$. We then have the following morphisms of DGAs:
$$
\H_{\dc}(A)(V) \la \z_{\dc}(A)(V) \to (A(V)),
$$
where $\z_{\dc}(A)(V)^n= \ker(\dc:A(V)^n \to A(V)^{n+1})$, with differential $d$, and $\H_{\dc}(A)(V)$ also has differential $d$. Since $\dc(a\cup b)= (\dc a)\cup b +(-1)^{\deg a} a \cup (db)$, these are indeed both DGAs. It follows from \cite{Simpson} Lemmas 2.1 and 2.2, using the $d\dc$ lemma instead of the $\pd\bar{\pd}$ lemma, that these morphisms are quasi-isomorphisms, and that $d=0$ on $\H_{\dc}(A(V))$.
\end{proof}
\end{theorem}

\begin{remark}
It follows from Remark \ref{toen} that this is equivalent to  \cite{KTP} Theorem 3.2.3, which states that the complex schematic homotopy type is formal.
\end{remark}

\begin{corollary}\label{quadratic}
The Lie algebra $\cL(\mypi)$ associated to $\mypi$ is quadratically presented (i.e. defined by equations of bracket length 2) as an element of $\widehat{\cN}(\redpi)$, and has a weight decomposition (as a pro-vector space), unique up to inner automorphism.
\begin{proof}
The functor $\mc_{\H(A)}$ is homogeneous, hence pro-representable, and \mbox{$\mc_{\H(A)} \to \ddef_{\H(A)}$} is \'etale, so $\mc_{\H(A)}$ is pro-represented by a hull for $\ddef_{\H(A)}$. By Theorem \ref{formal}, $\ddef_{\H(A)}$ is isomorphic to $\ddef_A$, which by Theorem \ref{vbeq} and Lemma \ref{repeq} is isomorphic to $R_{\rho_0}$. By Lemma \ref{gammahat}, this has hull $\mypi$. Therefore $\mypi$ pro-represents $\mc_{\H(A)}$, by the uniqueness of hulls.

Now, 
$$
\mc_{\H(A)}(\g) = \{ \omega \in \H^1(X, \mathpzc{g})\, |\, [\omega, \omega] =0 \in \H^2(X, \mathpzc{g})\}
$$
As in Remark \ref{toen}, we may replace $A$ by the $\redpi$-equivariant DGA $B:=A(O(\redpi))$. Letting $\bO$ denote the ind-local system on $X$ associated to the representation $O(\redpi)$, $H(A)$ then corresponds to the $\redpi$-equivariant DGA
$$
H(B)^n:=\H^n(X, \bO).
$$
The ind-local system $\bO$ was defined using the left $\redpi$-action on $O(\redpi)$, and the $\redpi$-action on $B$ is then defined using the right action.

Therefore $H(A)^n(V)=\H^n(X, \bO)\ten^{\redpi} V$, so
$$
\mc_{\H(A)}(\g) = \{ \omega \in \H^1(X, \bO)\ten^{\redpi} \g \,|\, [\omega, \omega] =0 \in \H^2(X, \bO)\ten^{\redpi} \g \}.
$$

Let 
$$
H_i:= \H^i(X, \bO)^{\vee} \in \widehat{\Rep}(\redpi).
$$
There are canonical isomorphisms 
$$
\Hom_{\widehat{\cN}(\redpi)}(L(H_1), \g) \cong \Hom_{\widehat{\Rep}(\redpi)}( \H_1, \g) \cong  \H^1(X, \mathpzc{g}),
$$
where $L$ denotes the free pro-nilpotent Lie algebra functor.

Now, the cup product
$$
\H^1(X,\bO)\ten \H^1(X,\bO) \xra{\cup} \H^1(X, \bO)
$$
gives a coproduct 
$$
\Delta: H_2 \to H_1 \hat{\ten} H_1 \to \hat{\bigwedge}^2 H_1 \subset  L(H_1).
$$

Finally, observe that $\mc_{\H(A)}(\g)$ is isomorphic to the set 
$$
\{\omega \in \Hom_{\widehat{\cN}(\redpi)}(L(H_1), \g) \,|\, \omega \circ \Delta(H_2)=0\}, 
$$
so
$$
\cL(\mypi) \cong L(H_1)/\Delta(H_2)
$$
is a quadratic presentation. 

If we set $H_1$ to have weight $-1$, and $H_2$ to have wight $-2$, then $\cL(\mypi)$ has a canonical weight decomposition arising from those on $H_1$ and $H_2$, since $\Delta$ preserves the weights. 
Note that a weight decomposition on a pro-finite-dimensional vector space is an infinite product, rather than an infinite direct sum. This decomposition is only unique up to inner automorphism, since the hull morphism is; this is equivalent to saying that we have not made a canonical choice of Levi decomposition.
\end{proof}
\end{corollary}

\begin{corollary}\label{criterion} Let $G$ be an arbitrary  reductive real algebraic group, acting on a real unipotent algebraic group $U$ defined by homogeneous equations, i.e. $\fu \cong \gr \fu$ as Lie algebras with $G$-actions. If  
$$
\rho_2: \pi_1(X,x) \to (U/[U,[U,U]]) \rtimes G
$$ 
is a Zariski-dense representation, then 
$$
\rho_1: \pi_1(X,x) \to (U/[U,U]) \rtimes G
$$
 lifts to a representation 
$$
\rho: \pi_1(X,x)\to U\rtimes G.
$$
\begin{proof}
 Observe that these representations correspond to $\pi_1(X,x)$-equivariant homomorphisms $\mypi \to U$. Let $\g$ be the Lie algebra associated to $\mypi$. We must show that the surjective map
$$
\rho_1:\g \to \fu/[\fu,\fu]
$$
lifts to $\fu$. Let $\g=L(V)/\langle W\rangle$, for $W \subset \wedge^2V$. Since $V \cong \gr_1\g=\g/[\g,\g]$, the morphism $\rho_1$ gives us a map $\theta:L(V) \to \gr \fu$, so it will suffice to show that $\theta$ annihilates $W$. But $\wedge^2\g  \to  \gr_2 \fu$ must send $u\wedge v$ to $[\rho_1(u),\rho_1(v)]=[\rho_2(u),\rho_2(v)]$, which annihilates $W$, as required.
\end{proof}
\end{corollary}

\begin{corollary}
For semisimple $\pi_1(X,x)$-representations $V_1, \ldots V_n$, with $n \ge 3$, the Massey products
$$
\H^1(\pi_1(X,x),V_1)\ten \H^1(\pi_1(X,x),V_2)\ten  \ldots \ten \H^1(V_n) \dashrightarrow \H^2(\pi_1(X,x), V_1\ten V_2 \ten  \ldots \ten V_n)
$$
are all zero.
\begin{proof} This follows from the observation that these maps all arise as quotients of higher obstruction maps for  quotients of $L(V_1\oplus V_2 \oplus \ldots \oplus V_n)$. Alternatively, it can be deduced directly from Theorem \ref{formal}, which implies that all the higher Massey products are zero on the cohomology of $X$ with semisimple coefficients.
\end{proof}
\end{corollary}

\begin{remarks}\label{compare} 
Note that Corollary \ref{quadratic} implies the results on the fundamental group  of \cite{DGMS}, of \cite{GM} and of \cite{malcev}. The pro-unipotent completion $\varpi_1(X,x)\ten \R$ studied in \cite{DGMS} is just the maximal quotient    of $\mypi$  on which $\pi_1(X,x)$ acts trivially. 

The problem considered in \cite{GM} (and generalised in \cite{Simpson}) is to fix a reductive representation $\rho_0:\pi_1(X,x) \to G(\R)$, and consider lifts $ \rho:\pi_1(X,x) \to G(A)$, for Artinian rings $A$. The hull of this functor is the functor
$$
A \mapsto \Hom_{\pi_1(X,x)}(\mypi, \exp(\g\ten \m_A)),
$$   
where $\g$ is the Lie algebra of $G$, regarded as the adjoint representation. It follows that this hull then has generators $\Hom_{\pi_1(X,x)}(\g, H_1)$, and relations 
$$
\Hom_{\pi_1(X,x)}(\g, H_2) \to S^2 \Hom_{\pi_1(X,x)}(\g, H_1)
$$ 
given by composing the coproduct and the Lie bracket.  

The statement of Corollary \ref{quadratic} is equivalent to saying that the relative Malcev completion of any Zariski-dense representation $\rho:\pi_1(X,x) \to G(\R)$, for $G$ reductive, is quadratically presented. In \cite{malcev} Theorem 13.14, this is proved  only for those $\rho$ which are polarised variations of Hodge structure, and no consequences are given. 
\end{remarks}

\begin{proposition}\label{crit}
If $\pi_1(X,x) = \Delta \rtimes \Lambda$, with $\Lambda$ acting reductively on the pro-unipotent completion $\Delta\ten \R$, then $\Delta\ten \R$ is quadratically presented.
\begin{proof} By Lemma \ref{semidirect}, we know that
$$
\rho_{\sharp}\mypi \cong (\Delta \ten \R)\by \Ru(\Lambda^{\alg}) \in \exp(\widehat{\cN}(\Lambda^{\red})),
$$
for $\rho:\redpi \to \Lambda^{\red}$. From Theorem \ref{quadratic}, we know that $\mypi$ is quadratically presented, hence so is $\rho_{\sharp}\mypi$. 

Now, to give a quadratic presentation for a Lie algebra $\g$ is equivalent to giving a homomorphism $\theta:\gr(\g) \to \g$ from its associated graded Lie algebra such that \mbox{$\gr (\theta):\gr (\g) \to \gr (\g)$} is   the identity,  provided that $\gr (\g)$ is quadratic as a graded Lie algebra. Now, if $\g \oplus \fh$ is quadratically presented, then so is $\gr(\g \oplus \fh)=\gr(\g)\oplus \gr(\fh)$, hence so is $\gr(\g)$. Taking the composition
$$
\gr(\g) \to \gr(\g)\oplus \gr(\fh) \xra{\theta} \g \oplus \fh \to \g
$$
then gives a quadratic presentation for $\g$.

Combining these results, we see that $\Delta\ten \R$ must be quadratically presented.
\end{proof}
\end{proposition}

\begin{remark}\label{hain}
Since \cite{Simpson} Lemma 4.5 states that properly rigid reductive representations underlie variations of Hodge structure, Proposition \ref{crit} can be deduced directly from \cite{malcev} whenever the composition $\pi_1(X,x) \to \L \to \Aut(\Delta\ten \R)$  is  properly rigid and reductive.
\end{remark}

\begin{example} Let $\fh = \R^2 \oplus \R$,  with Lie bracket $[\fh,\R]=0$ and $[u,v]=u \wedge v \in \R$ for $u,v \in \R^2$, so $\exp(\fh)$ is isomorphic to the real three-dimensional Heisenberg group. The Campbell-Baker-Hausdorff formula enables us to regard $\exp(\fh)$ as the group with underlying set $\fh$ and product
$
a\cdot b=a + b +\half [a,b],
$
since all higher brackets vanish. It then follows that the lattice
$$
H:=\exp(\Z^2 \oplus  \half \Z),
$$
is closed under this multiplication, so forms a discrete group, with $H\ten \R = \exp(\fh)$. Now, $\SL_2(\Z)$ acts on $H$ by the formula:
$$
A(v,w) := (Av, (\det A) w)=(Av,w),
$$
for $v \in \Z^2$ and $w \in \half \Z$. 

Let $\L:=\Z^2=\Z a\oplus \Z b$ act on $H$ via the homomorphism $\vartheta:\L \to \SL_2(\Z)$ given by 
$$
\vartheta(a)=\vartheta(b)=M:= \begin{pmatrix} 2 & 3 \\ 1 &2 \end{pmatrix}.
$$
Then the Zariski closure of this representation is isomorphic to $\bG_m(\R)=\R^*$, being the torus in $\SL_2(\R)$ containing $M$, so the action of $\L$ on $\fh$ is reductive. Since $\fh$ is not quadratically presented, the group 
$\Gamma:=H \rtimes \SL_2(\Z)$ 
cannot be the fundamental group of any compact K\"ahler manifold.

Note that \cite{DGMS} cannot be used to exclude this group: since  the commutator $[a,(v,0)]= (Mv-v,0)$ and $M-I$ is non-singular, $[\Gamma,\Delta]$ is of finite index in $\Delta$, so $\Gamma\ten \R=\L\ten \R$, which is quadratically presented.

Furthermore, this result cannot be obtained by substituting Remark \ref{hain} and \cite{malcev} Theorem 13.14 for Corollary \ref{quadratic}, since $\vartheta$ is not rigid.    

Alternatively, we could use Corollary \ref{criterion} to prove that $\Gamma$ is not a K\"ahler group. Let $G=\bG_m(\R)$, $\fu=L(\R^2)$ and $U=\exp(\fu)$. Observe that $\fh \cong \fu/[\fu,[\fu,\fu]]$, and let
$$
\rho_2: H \rtimes \L \to \exp(\fh) \rtimes \bG_m(\R),
$$
be given by combining the standard embedding with $\vartheta$.

Since all triple commutators vanish in $H$, this does not lift to a representation
$$
\rho: H \rtimes \L \to U \rtimes \bG_m(\R).
$$
\end{example}

%\newpage
\bibliographystyle{alpha}
\addcontentsline{toc}{section}{Bibliography}
\bibliography{references.bib}
\end{document}